\documentclass[final, 12pt]{article}
\usepackage{amsfonts}
\usepackage{makeidx}
\usepackage{amsmath}
\usepackage{amssymb}
\usepackage{bm}
\usepackage{cite}
\usepackage{showkeys}
\usepackage[english]{babel}
\usepackage{dblaccnt}
\usepackage{accents}
\usepackage{graphicx}
\usepackage{psfrag}
\usepackage{subfig}
\usepackage{color}
\usepackage{float}
\usepackage{amscd}
\usepackage[mathscr]{euscript}
\usepackage{lipsum}
\usepackage{epstopdf}
\usepackage{algorithm}
\usepackage{algpseudocode}

\newtheorem{defi}{Definition}

\newtheorem{remark}[defi]{Remark}
\begin{document}

\title{Convexification  numerical algorithm for a 2D inverse
scattering problem with backscatter data}
\author{Trung Truong\thanks{
Department of Mathematics, Kansas State University, USA; (\texttt{%
trungt@ksu.edu})} \and Dinh-Liem Nguyen\thanks{%
Department of Mathematics, Kansas State University, USA; (\texttt{%
dlnguyen@ksu.edu}) } \and Michael V. Klibanov\thanks{
Department of Mathematics and Statistics, University of North Carolina at
Charlotte, USA; (\texttt{mklibanv@uncc.edu}) } }
\date{}
\maketitle

\begin{abstract}
This paper is concerned with the inverse scattering problem which aims to
determine the spatially distributed dielectric constant coefficient of the
2D Helmholtz equation from multifrequency backscatter data associated with a single
direction of the incident plane wave. We propose a globally convergent
convexification numerical algorithm to solve this nonlinear and ill-posed
inverse problem. The key advantage of our method over conventional
optimization approaches is that it does not require a good first guess about
the solution. First, we eliminate the coefficient from the Helmholtz
equation using a change of variables. Next, using a truncated expansion with
respect to a special Fourier basis, we approximately reformulate the inverse
problem as a system of quasilinear elliptic PDEs, which can be numerically
solved by a weighted quasi-reversibility approach. The cost functional for
the weighted quasi-reversibility method is constructed as a Tikhonov-like
functional that involves a Carleman Weight Function. Our numerical study
shows that, using a version of the gradient descent method, one can find the
minimizer of this Tikhonov-like functional without any advanced \emph{a priori} knowledge
about it.
\end{abstract}

\sloppy

\textbf{Keywords. } inverse scattering, numerical reconstruction,
convexification, backscatter data, Carleman weight function, coefficient
identification

\bigskip

\textbf{AMS subject classification. } 35R30, 78A46, 65C20

\section{Introduction}

\label{sec1}

%

Consider the scattering problem for a penetrable inhomogeneous medium in $%
\mathbb{R}^{2}$. Below ${x}=(x_{1},x_{2})^{T}\in \mathbb{R}^{2}.$ We assume
that the scattering object, which occupies a bounded domain in $\mathbb{R}%
^{2}$, is characterized by the spatially distributed dielectric
constant $\varepsilon_{r}( x) =1+a( x)$,
where the function $a(x)$ has a compact support. In this paper we
are particularly interested in the case of $a(x)\geq 0$ that typically
appears in applications of non-destructive testing and explosive detection,
see for instance~\cite{Bukhg1981, Klib95, Klib97} for a
similar assumption. Suppose that the object is illuminated by the downward propagating incident
plane wave $u_{\mathrm{in}}(x,k)=\exp(ik(d_1x_1+d_2x_2))$, where $d_1^2+d_2^2 = 1,  d_2<0$, the propagation direction
$(d_1,d_2)^T$ is fixed, and $k$ is the wavenumber. Then there arises the scattered
wave, and the total wave $u(x,k)$ which is the sum of the incident wave and the
scattered wave is governed by the Helmholtz equation as 
\begin{align}
\label{Helm}
& \Delta u+k^{2}(1+a(x))u=0,\quad x\in \mathbb{R}^{2},   \\
\label{radiation}
& \lim_{|x|\rightarrow \infty }|x|\left( \frac{\partial (u-u_{\mathrm{in}})}{%
\partial |x|}-ik(u-u_{\mathrm{in}})\right) =0.
\end{align}%
The scattered wave $u-u_{\mathrm{in}}$
satisfies the Sommerfeld radiation condition~\eqref{radiation}, which
guarantees that it behaves like a spherically outgoing wave far away from
the scattering object. It is well known that the scattering 
problem~\eqref{Helm}--\eqref{radiation} has a unique solution $u$, see~\cite{Colto2013}. Now let $R>0$ and consider 
\begin{equation*}
\Omega =(-R,R)^{2},\quad \Gamma =(-R,R)\times \{R\}.
\end{equation*}%
Assume that the scatterer as well as the support of the coefficient $a(x)$
are contained in $\Omega $, and that these objects do not intersect with $%
\partial \Omega $. Let $\underline{k}$ and $\overline{k}$ be positive
constants such that $\underline{k}<\overline{k}$. We consider the following
inverse problem.

\begin{figure}[h]
\centering
\hspace{-2.5cm} \includegraphics[width=6cm]{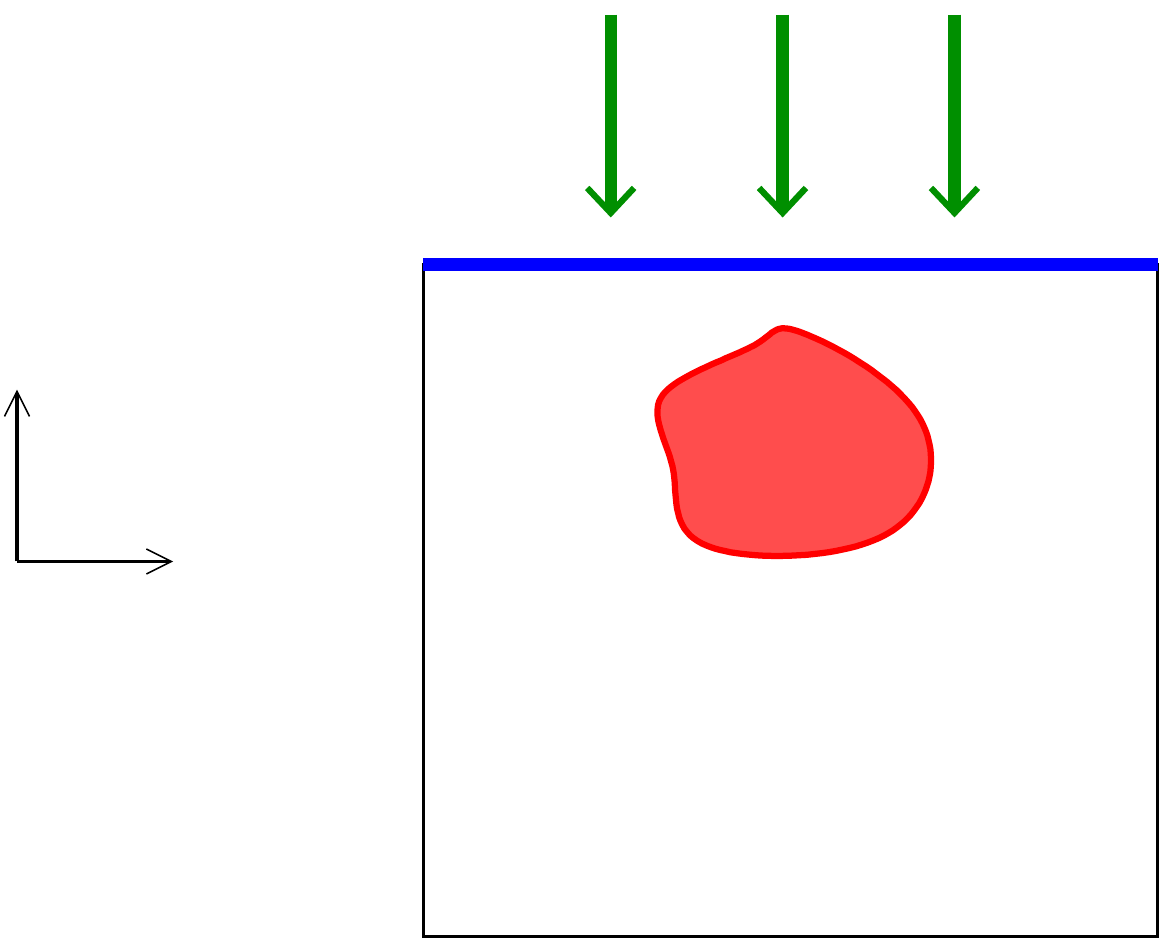}
\caption{Schematic of the inverse scattering from a penetrable bounded
object characterized by the function $a(x)$. The incident plane wave
propagates downward toward the scattering object. The backscatter data are
measured on the top boundary $\Gamma$ of the computational domain $\Omega$.}
\label{fg0}
\end{figure}

\textbf{Inverse Problem.} Assume that we are given the multi-frequency
backscatter Cauchy data 
\begin{align}
g_{0}(x,k)& :=u(x,k),\quad \text{for }{x}\in \Gamma ,k\in \lbrack \underline{%
k},\overline{k}],  \label{data2} \\
g_{1}(x,k)& :=\frac{\partial u}{\partial x_{2}}(x,k),\quad \text{for }{x}\in
\Gamma ,k\in \lbrack \underline{k},\overline{k}],  \label{data3}
\end{align}%
where the total wave $u(x,k)$ is generated by incident plane waves with a fixed
propagation direction. 
Determine the function $a(x)$ in~\eqref{Helm} for $x\in
\Omega $, see also Figure~\ref{fg0} for a schematic diagram of the
measurement arrangement in the inverse problem.

Uniqueness theorem for this inverse problem can be currently proven
only in the case when the right hand side of equation (\ref{Helm}) does
not equal to zero in $\overline{\Omega}$. This can be done by the
so-called Bukhgeim-Klibanov method, which was originated in \cite{Bukhg1981}
and is based on applications of Carleman estimates to coefficient inverse
problems, see, e.g.~\cite{Beili2012,KT,Ksurvey} for
this method. In addition, uniqueness of the
approximate problem can be proven when the truncated Fourier series
for~\eqref{Helm} is used for that approximation, see, e.g. Theorem 3.2 in~\cite{Khoa2019}.

\begin{picture}(10,3)
\put(245,198){\small {\color{blue}$\Gamma$}}
\put(198,123){\small $\Omega$}
\put(263,138){\small $a = 0$}
\put(191,163){\small $a(x)$}
\put(195,233){\small {\color{black} $u_{\mathrm{in}}$}}
\put(113,151){\small $x_1$}
\put(90,176){\small $x_2$}
\end{picture}

This inverse problem belongs to a wider class of coefficient inverse
scattering problems which in general aim to recover information about the
coefficient $a(x)$ (e.g. its support and/or its values) from the knowledge
of the scattered wave generated by a number of incident waves. Inverse
scattering problems occur in many applications, including non-destructive
testing, explosive detection, medical imaging, radar imaging and geophysical
exploration. There is a vast literature about theoretical results and
numerical solution to inverse scattering problems, see for instance~\cite%
{Colto2013} and references therein. Due to the interest of this paper, we
discuss only some numerical methods. The conventional approach is based on
the optimization based methods, see, e.g.~\cite{Bakus2004,
Chavent,Gonch1,Gonch2,Engl1996}. However, it is well known that these
methods may suffer from multiple local minima and ravines and their
convergence analysis is also unknown in many situations. An important
attempt in overcoming the drawbacks of the optimization based methods is the
qualitative approach which aims to compute the geometry of the scattering
object or the support of the coefficient $a(x)$. We refer to~\cite{Cakon2006, Kirsc2008, Colto2013} and references therein for the development
of qualitative methods in solving inverse scattering problems. Although one
may be able to avoid local minima or the use of advanced \emph{a priori}
information of the solution, still only geometrical information of the
scatterer can be reconstructed with qualitative methods. Furthermore, these
methods typically require muti-static data which are sometimes not available
in practical applications.

The numerical method proposed in this paper is an extended study from a
recent new approach called globally convergent numerical methods (GCNM) for
solving coefficient inverse problems. We say that a numerical method
for a nonlinear ill-posed problem converges \emph{globally} if there is a
rigorous guarantee that it delivers points in a sufficiently small
neighborhood of the exact solution of this problem without any advanced
knowledge of this neighborhood. The GCNM typically aims to reconstruct a
coefficient in an inverse scattering problem using scattering data either
for a single direction of the incident plane wave, or, most recently, for
many locations of the point source but at a fixed single frequency \cite%
{Khoa2019}. An interesting feature of GCNM is that in all cases the
data are non over-determined. The latter means that the number $m$
of free variables in the data equals the number $n$ of free
variables in the unknown coefficient, $m=n$. The main advantage of
the GCNM is that any version of it avoids the local minimum problem suffered
by optimization based methods. Still, any version of GCNM holds the above
indicated global convergence property. We refer to~\cite{Beili2012,
Kliba2018, Nguye2017, Koles2017, Nguye2018} and references therein for
theoretical results as well as numerical and experimental data study of the
first type of GCNM. 

The method of this paper is inspired by the second type of the GCNM,
which is called convexification. The development of the
convexification has started in 1995 and 1997 by Klibanov \cite{Klib95,Klib97}
and continued since then in \cite{BK,KT,KK}. However, those were mostly
analytical works since some obstacles existed at that time on the path to
the numerical implementation, although see some numerical results 
for the one-dimensional case in \cite{KT}. 
Fortunately, in 2017 the work \cite{Bak} has eliminated those
obstacles. This generated a number of more recent publications on the
convexification~\cite{Kliba2018a, Kliba2019a, Kliba2017, Kliba2019c, Kliba2019d, Khoa2019}, which contain both a rigorous convergence
analysis and numerical results. In particular, publications \cite{Kliba2018a, Kliba2019a} are
about the verification of the convexification on experimental data. 

The central idea of the convexification is to construct of a
globally  convex weighted Tikhonov-like functional with the Carleman
Weight Function in it. The idea of the use of the Carleman Weight Function
is an  unexpected consequence of the original idea of the
Bukhgeim-Klibanov method~\cite{Bukhg1981}, which was originally aimed
only for proofs of uniqueness theorems for coefficient inverse problems. The
final step of the convergence analysis of the convexification consists in
the proof of the global convergence of the gradient projection method to the
exact solution, as long as the level of noise in the data tends to zero. 
We also refer to another version of the convexification, which has
started in the work~\cite{Baud1} and has been continued in \cite{Baud2,BBS,LN}. 
Carleman Weight Functions are also a crucial element of
these works. The main difference between these publications and our method
is that it is assumed in \cite{Baud1,Baud2,BBS,LN} that the initial
condition in a hyperbolic/parabolic PDE is not vanishing in $\overline{\Omega}$, 
which unlike our case of the zero right hand side of
equation~\eqref{Helm}. 

As to this present paper, our first step is to eliminate the
coefficient from the Helmholtz equation using a change of variables. Next,
using a truncated Fourier expansion for a function generated by the total
wave field, we approximately reformulate the inverse problem as the Cauchy
problem for a system of quasilinear elliptic PDEs. The Cauchy boundary data
are as follows: on a part of the boundary both Dirichlet and Neumann
boundary data are given and no data are given on the rest of the boundary.
We then propose a weighted quasi-reversibility method to solve the 
problem. Inspired by the concept of the convexification, the cost functional
in that weighted quasi-reversibility method contains a Carleman Weight
Function. This function plays the decisive role in the numerical performance
of the method. A method of gradient descent type is explored to find the
global minimizer of the cost functional without using any advanced \emph{a priori}
information about it.

Comparing with the above cited recent works on the convexification, the new
features of this work are that firstly our algorithm exploits the new
Fourier basis in~\cite{Kliba2017} to solve a multi-dimensional inverse
problem for the Helmholtz equation with multifrequency data and a
single direction of the incident plane wave. The latter is mostly related
to~\cite{Kliba2018a} in which, however, only the one-dimensional
version of the inverse problem has been studied. Using the new Fourier basis
from~\cite{Kliba2017}, the 3D inverse problem for the Helmholtz equation
with data generated by a moving source (at a fixed frequency) has been also
studied in~\cite{Khoa2019}. Secondly, a modification during the iteration of
the gradient descent method is applied to help the cost functional converge
faster. More precisely, we solve the direct problem to update some functions 
during the iterations of the gradient descent method. Thirdly,
unlike the previous works~\cite{Kliba2018a, Kliba2019a, Kliba2017, Kliba2019c, Kliba2019d, Khoa2019}, 
the reconstruction
algorithm proposed in this paper works without using any data completion
process, and the numerical study covers challenging cases of scattering
objects of different shapes which are characterized by different 
values of the dielectric constant. We also want to mention that the
implementation of the method uses a full $H^{2}$ term instead of $L^{2}$ or $%
H^{1}$ terms as in the previous works cited above and does not need any
cut-off and averaging procedures during the iteration in the algorithm. 

The convergence analysis of the method of this paper will be addressed in an
incoming publication. To be more precise, we now roughly (i.e.
without some details) specify what kind of theorems will be proven in that
publication. Analogs of these theorems for the one-dimensional case can be found in~\cite{Kliba2018a}.
Those theorems claim:

\begin{enumerate}
\item  The strict convexity of the weighted Tikhonov-like functional 
$\mathcal{J}(W)$ in (\ref{cost}) on the ball $B( M)
\subset \mathcal{X}$ of the radius $M$, see (\ref{200}) and (\ref{201}).
 The strict convexity will be proven for sufficiently large values of the
parameter $\lambda \geq \lambda ( M) \geq 1$ in the
Carleman Weight Function (\ref{cw}).

\item Existence and uniqueness of the minimizer of the functional $%
\mathcal{J}(W)$ on $\overline{B( M) }$ for $%
\lambda \geq \lambda ( M)$.

\item Convergence of the gradient projection method of the
minimization of the functional $\mathcal{J}(W)$ on $\overline{B( M)}$
 to the exact solution of that approximate
coefficient inverse problem if starting from an arbitrary point of
$B( M)$. That convergence takes place as long as the
level of noise in the data tends to zero. 
\end{enumerate}

Since the radius $M>0$ of the ball $B( M)$
is an arbitrary number, then this is the desired global
convergence property, as defined above. Note that even though the theory
requires the parameter $\lambda$ to be sufficiently large, the
numerical experience of this and all previous publications about the
convexification~\cite{Kliba2018a, Kliba2019a, Kliba2017, Kliba2019c, Kliba2019d, Khoa2019} 
shows that actually reasonable values of 
$\lambda$ provide accurate solutions of considered inverse problems.

The paper is structured as follows. The second section is dedicated to the
formulation of the inverse problem as an approximate quasilinear elliptic
PDE system. The numerical reconstruction method for solving the inverse
problem is proposed in Section 3. The implementation and numerical examples
of the reconstruction method are presented in Section 4. Finally, Section 5
contains a summary discussion of this work.

\section{An approximate elliptic PDE formulation}

In this section we reduce our inverse problem to the Cauchy problem
for a system of quasilinear elliptic PDEs that we will be studying using a
quasi-reversibility approach in the next section. The main ideas for
deriving the formulation are using truncated Fourier expansion in $L^{2}(%
\underline{k},\overline{k})$ and eliminating the coefficient $a(x)$ from the
scattering problem. Setting $k_{0}=(\underline{k}+\overline{k})/2$ we first
need the following important Fourier basis of $L^{2}(\underline{k},\overline{%
k})$ that was introduced in~\cite{Kliba2017} 
\begin{equation*}
\psi _{n}(k)=(k-k_{0})^{n-1}e^{k-k_{0}},\quad k\in (\underline{k},\overline{k%
}),\quad n=1,2,\dots 
\end{equation*}%
Applying the Gram--Schmidt process to $(\psi _{n})$ we obtain an orthonormal
basis $\left\{ \Phi _{n}\right\} _{n=1}^{\infty },$ which has the following
properties, also, see~\cite{Kliba2017}:\newline
i) $\Phi _{n}\in C^{\infty }[\underline{k},\overline{k}]$ for all $n=1,2,...$
\newline
ii) The matrix $D=[d_{mn}]$, where $m,n=1,\dots ,N$ and 
\begin{equation*}
d_{mn}=\int_{\underline{k}}^{\overline{k}}\Phi _{n}^{\prime }(k)\Phi
_{n}(k)dk,
\end{equation*}%
is invertible with $d_{mn}=1$ for $m=n$ and $d_{mn}=0$ for $m>n$.

Now setting 
\begin{equation}
p(x,k)=\frac{u(x,k)}{u_{\mathrm{in}}(x,k)},  \label{p}
\end{equation}%
and substituting in~\eqref{Helm} we obtain 
\begin{equation}
\Delta p(x,k)+k^{2}a(x)p(x,k)-2ik\partial _{x_{2}}p(x,k)=0.  \label{hew}
\end{equation}%
Now suppose that $p(x,k)$ is nonzero for all $x\in \Omega ,k\in \lbrack 
\underline{k},\overline{k}]$. We define $v(x,k)$ as 
\begin{equation}
v(x,k)=\frac{\log (p(x,k))}{k^{2}},  \label{v}
\end{equation}%
where $\log $ is the principal logarithm. We also assume that $v(x,k)$ is
continuous and differentiable for all $x\in \Omega ,k\in \lbrack \underline{k%
},\overline{k}]$. We refer to~\cite{Khoa2019,Kliba2019a,Kliba2019d} for the
definition of the complex logarithm for a similar change of variables using
a high frequency asymptotic behavior for the total field in $\mathbb{R}^{3}$. 
Next, this definition was extended in \cite{Khoa2019} to non high
values of $k$ as long as $v(x,k)\neq 0$ for those
values. To what we know, that asymptotic behavior is not established yet
for the two-dimensional case. At the same time, in our numerical
studies, we do not see any discontinuity problem with the principal $\log$.

Using~\eqref{v} we substitute $p=\exp (k^{2}v)$ in~\eqref{hew} and rewrite (%
\ref{hew}) in terms of $v(x,k)$ as follows 
\begin{equation}
\Delta v(x,k)+k^{2}\nabla v(x,k)\cdot \nabla v(x,k)-2ik\partial
_{x_{2}}v(x,k)+a(x)=0.  \label{hev}
\end{equation}%
We now eliminate $a(x)$ by differentiating (\ref{hev}) with respect to $k$ 
\begin{equation}
\Delta (\partial _{k}v)+2k\nabla v\cdot \nabla (v+k\partial _{k}v)-2i\left(
\partial _{x_{2}}v+k\partial _{x_{2}}\partial _{k}v\right) =0.  \label{hedk}
\end{equation}%
Let $N\in \mathbb{N}$ be sufficiently large. We approximate the function $%
v(x,k)$ in~\eqref{v} and its partial derivative $\partial _{k}v(x,k)$ using
the truncated Fourier series as 
\begin{equation}
v(x,k)=\sum_{n=1}^{N}v_{n}(x)\Phi _{n}(k),\quad \partial
_{k}v(x,k)=\sum_{n=1}^{N}v_{n}(x)\Phi _{n}^{\prime }(k),  \label{100}
\end{equation}%
where the coefficients $v_{n}(x)$ are given by 
\begin{equation}
v_{n}(x)=\int_{\underline{k}}^{\overline{k}}v(x,k)\Phi _{n}(k)dk.
\label{coef}
\end{equation}%
Using two truncated series (\ref{100}), we approximate (\ref{hedk}) by 
\begin{align*}
 \sum_{n=1}^{N}\Phi _{n}^{\prime }(k)\Delta
v_{n}(x)+2k\sum_{n=1}^{N}\sum_{l=1}^{N}\Phi _{n}(k)(\Phi _{l}(k)+k\Phi
_{l}^{\prime }(k))\nabla v_{n}(x)\cdot \nabla v_{l}(x) \\
-2i\sum_{n=1}^{N}(\Phi _{n}(k)+k\Phi _{n}^{\prime }(k))\partial
_{x_{2}}v_{n}(x)& =0.
\end{align*}%
For each $m=1,\dots ,N$, multiplying both sides of the above equation by $%
\Phi _{m}(k)$ and integrating with respect to $k$ over $[\underline{k},%
\overline{k}]$, we obtain 
\begin{multline}
\sum_{n=1}^{N}\left( \int_{\underline{k}}^{\overline{k}}\Phi _{m}(k)\Phi
_{n}^{\prime }(k)dk\right) \Delta v_{n}(x)  \label{appr} \\
+\sum_{n=1}^{N}\sum_{l=1}^{N}\left( 2k\int_{\underline{k}}^{\overline{k}%
}\Phi _{m}(k)\Phi _{n}(k)\left[ \Phi _{l}(k)+k\Phi _{l}^{\prime }(k)\right]
dk\right) \nabla v_{n}(x)\cdot \nabla v_{l}(x) \\
-\sum_{n=1}^{N}\left( 2i\int_{\underline{k}}^{\overline{k}}\Phi _{m}(k)\left[
\Phi _{n}(k)+k\Phi _{n}^{\prime }(k)\right] dk\right) \partial
_{x_{2}}v_{n}(x)=0.
\end{multline}%
Considering two $N\times N$ matrices defined as 
\begin{align*}
D& =(d_{mn}),\quad d_{mn}=\int_{\underline{k}}^{\overline{k}}\Phi
_{m}(k)\Phi _{n}^{\prime }(k)dk, \\
S& =(s_{mn}),\quad s_{mn}=-2i\int_{\underline{k}}^{\overline{k}}\Phi _{m}(k)%
\left[ \Phi _{n}(k)+k\Phi _{n}^{\prime }(k)\right] dk,
\end{align*}%
and an $N\times N$ block matrix $B=(B_{mn})$, each block $%
B_{mn}=(b_{mn}^{(l)})_{l}$ is an $N\times 1$ matrix defined as 
\begin{equation*}
b_{mn}^{(l)}=2k\int_{\underline{k}}^{\overline{k}}\Phi _{m}(k)\Phi _{n}(k)%
\left[ \Phi _{l}(k)+k\Phi _{l}^{\prime }(k)\right] dk,
\end{equation*}%
we can rewrite (\ref{appr}) as a system of PDEs for the vector valued
function $V(x)=[v_{1}(x)\ v_{2}(x)\ \dots \ v_{N}(x)]^{T}$ 
\begin{equation}
D\Delta V(x)+B\partial _{x_{1}}V(x)\bullet \partial _{x_{1}}V(x)+B\partial
_{x_{2}}V(x)\bullet \partial _{x_{2}}V(x)+S\partial _{x_{2}}V(x)=0.
\label{appr2}
\end{equation}%
Here the operator $\bullet $ is defined as follows: If $P=(P_{m})$ is an $%
N\times 1$ block matrix, each block $P_{m}$ is an $N$-dimensional column
vector and $V$ is an $N$-dimensional column vector then $P\bullet V$ is an $N
$-dimensional column vector given by 
\begin{equation*}
P\bullet V=%
\begin{bmatrix}
P_{1}\cdot V \\ 
P_{2}\cdot V \\ 
\vdots  \\ 
P_{N}\cdot V%
\end{bmatrix}%
.
\end{equation*}%
Defining 
\begin{equation*}
\mathcal{Q}(V)=D\Delta V+B\partial _{x_{1}}V\bullet \partial
_{x_{1}}V+B\partial _{x_{2}}V\bullet \partial _{x_{2}}V+S\partial _{x_{2}}V,
\end{equation*}%
we are able to approximately reformulate the inverse problem as the
Cauchy problem for the following system of quasilinear elliptic PDEs:%
\begin{align}
\mathcal{Q}(V)& =0\quad \text{in}\ \Omega ,  \label{eqn} \\
V& =G_{0}\quad \text{on}\ \Gamma ,  \label{dd} \\
\partial _{x_{2}}V& =G_{1}\quad \text{on}\ \Gamma ,  \label{nd}
\end{align}%
where $G_{0}$ and $G_{1}$ can be computed from the given boundary data $g_{0}
$ and $g_{1}$ in (\ref{data2})--(\ref{data3}) using~\eqref{p}, \eqref{v} and~%
\eqref{coef}. If we can find $V$ by solving problem~\eqref{eqn}--\eqref{nd},
the coefficient of interest $a(x)$ can be approximately recovered from~%
\eqref{hev}.

\begin{remark}
We emphasize that the reconstruction algorithm we study in the next section
for solving problem~\eqref{eqn}--\eqref{nd} only needs the backscatter data
on $\Gamma $. In contrast, the convexification method of above cited 
papers~\cite{Kliba2018a, Kliba2019a, Kliba2017, Kliba2019c, Kliba2019d, Khoa2019}, 
 one has to artificially complete the backscatter
data on the other boundaries of $\Omega $ for a better stability of
computations. On the other hand, the forthcoming analytical
results that are mentioned in Introduction for this paper are valid with the Carleman Weight Function~\eqref{cw}
 only if the Dirichlet data for the system (\ref{eqn}) are known on the
entire boundary $\partial \Omega $ rather than just on its part $
\Gamma$, i.e. they are valid for those completed data. Thus, our
claim in the first sentence of this Remark is based only on our numerical
observation and is not supported by the theory. Nevertheless, this numerical
observation  emphasizes the stability property of our method.
\end{remark}

\begin{remark}
It is well known that the Cauchy problem for an elliptic equation is
unstable. Thus, we actually construct a regularization method of solving
this problem for our case. A similar numerical method was constructed in 
\cite{Klib2015} for ill-posed Cauchy problems for a wide class of single
quasilinear PDEs, including the elliptic one. However, the Carleman Weight
Function used in \cite{Klib2015}  for the elliptic case is inconvenient for
the numerical implementation since it depends on two large parameters,
instead of just one in our case of (\ref{cw}).
\end{remark}

\section{A numerical reconstruction algorithm}

We solve problem~\eqref{eqn}--\eqref{nd} using the weighted
quasi-reversibility method. We first make a change of variables to have
homogeneous boundary conditions on $\Gamma $. Let $F$ be a vector valued
function which satisfies the boundary conditions (\ref{dd})--(\ref{nd}). We
call $F$ the data carrier and its construction is detailed in the numerical
study section. Assuming $V$ is the solution of problem~\eqref{eqn}--%
\eqref{nd}, we define 
\begin{equation*}
W=V-F.
\end{equation*}%
Then $W$ satisfies the homogeneous boundary conditions on $\Gamma $, that
is, 
\begin{equation*}
W=\partial _{x_{2}}W=0\quad \text{on}\ \Gamma .
\end{equation*}%
Define the function space $\mathcal{X}$ as 
\begin{equation}
\mathcal{X}=\left\{ W\in \lbrack H^{2}(\Omega )]^{N},\ W=\partial
_{x_{2}}W=0\ \text{on}\ \Gamma \right\}   \label{200}
\end{equation}%
with its associated norm 
\begin{equation*}
\lVert W\rVert _{\mathcal{X}}=\left( \sum_{n=1}^{N}\lVert w_{n}\rVert
_{H^{2}(\Omega )}^{2}\right) ^{\frac{1}{2}},\quad \text{where }%
W(x)=[w_{1}(x)\ w_{2}(x)\ \dots \ w_{N}(x)]^{T}.
\end{equation*}

%

Let $M>0$ be an arbitrary number. Define the ball $%
B( M) \subset \mathcal{X}$ as
\begin{equation}
B( M) =\left\{ W\in \mathcal{X}:\left\Vert W\right\Vert _{%
\mathcal{X}}<M\right\} \subset \mathcal{X}.  \label{201}
\end{equation}%
Next, we define the weighted Tikhonov-like functional $\mathcal{J}
:B( M) \rightarrow \mathbb{R}$ as
\begin{equation}
\mathcal{J}(W)=\int_{\Omega }\lvert \mathcal{Q}(W+F)\rvert ^{2}\varphi
^{2}dx+\rho \lVert W\rVert _{\mathcal{X}}^{2}+\alpha _{1}\int_{\Gamma
}\lvert W\rvert ^{2}dx+\alpha _{2}\int_{\Gamma }\left\vert \partial _{x_{2}}{%
W}\right\vert ^{2}dx,  \label{cost}
\end{equation}%
where 
\begin{equation}
\varphi (x)=e^{-\lambda (x_{2}-s)^{2}}  \label{cw}
\end{equation}%
is a Carleman Weight function and $\lambda \geq 1$ and $s>R$ are constants.
From our numerical experience, the regularization terms involving $\alpha_1$ and $\alpha_2$
help us obtain  better stability for the computation although they are not needed 
in the  theory of convexification methods in previous studies~\cite{Kliba2018a, Kliba2019a, Kliba2017, Kliba2019c, Kliba2019d, Khoa2019}.
Below we focus on the minimization of the functional $\mathcal{J}(W)$ on the
ball $B( M) \subset \mathcal{X}$ defined in (\ref{201}). 
As stated in Introduction, the use of the Carleman Weight Function is
inspired by convexification methods whose different versions are described
in the above cited publications. For the Carleman Weight Function $\varphi $
in~\eqref{cw}, a Carleman estimate for the Laplacian has been proved in~\cite{Kliba2019c}, 
where the Dirichlet boundary condition is given on the entire
boundary of $\Omega $, which requires an artificial complement of
the backscatter data given only on the part $\Gamma $ of $
\partial \Omega$. Actually, assigning the Dirichlet data on the
entire boundary $\partial \Omega$, one provides an additional
stability property to the method. Recall that (see Remark 1) it is our
numerical observation that our algorithm only requires the backscatter data
on the top boundary $\Gamma$ of $\Omega $.

An interesting numerical observation is that our algorithm converges
and provides better reconstruction results with the Carleman Weight
Function, as compared with the case when this function is absent in 
(\ref{cost}), i.e. when $\lambda =0$ in (\ref{cw}): compare 
Figure~\ref{fi1} and Figure~\ref{fi0} for numerical results with and without the
Carleman Weight Function. The two parameters $\lambda $ and $s$ along with
the regularization parameters $\rho$, $\alpha_{1}$, $\alpha_{2}$ will be
chosen numerically in the implementation of the algorithm. We find the
solution $W$ as the global minimizer of $\mathcal{J}(W)$ using a method of
gradient descent type. We point out that even though we can prove
the global convergence on $B(M) $ of the gradient
projection method rather than of the gradient descent method, still our
numerical observation is that the latter method has good convergence
properties. The same observation took place in all previous publications
about the convexification where numerical results were presented~\cite{Kliba2018a, Kliba2019a, Kliba2017, Kliba2019c, Kliba2019d, Khoa2019}. 
This is a quite useful observation since the numerical
implementation of the gradient descent method is much simpler than the one
of the gradient projection method.

In the following we describe our numerical algorithm for finding the
coefficient $a(x)$ for the inverse problem in~\eqref{data2}--\eqref{data3},
in which finding the minimizer $W$ of $\mathcal{J}(W)$ is one of the main
components of the algorithm.

\begin{remark}
In this algorithm, we recall that the capital letter notations, for example $%
W(x)=[w_{1}(x)\ w_{2}(x)\ \dots \ w_{N}(x)]^{T}$, are vector valued
functions whose components are Fourier coefficients of the corresponding
scalar function 
\begin{equation*}
w(x,k)=\sum_{n=1}^{N}w_{n}(x)\Phi _{n}(k)
\end{equation*}%
with the normal letter notations. We also prescribe a tolerance which forces
our iteration to stop after the cost functional no longer decreases much.
The tolerance will be chosen numerically in the implementation of the
algorithm.
\end{remark}

\begin{centering}
\bf{The numerical reconstruction algorithm}
\end{centering}

\paragraph{Step 1.}

Construct the data carrier $F$. Set the initial guess $V_{0}:=F$ then
proceed into the main iteration (\textbf{Step 2}).

%
%

\paragraph{Step 2a.}

Set $W_n := V_n - F$, compute the cost functional $\mathcal{J}(W_n)$ and its
gradient $\nabla\mathcal{J}(W_n)$.

\begin{itemize}

\item If $n\geq1$ and $|\mathcal{J}(W_n) - \mathcal{J}(W_{n-1})| < tolerance$%
, stop the iteration and move to \textbf{Step 3}.

\item Otherwise proceed to \textbf{Step 2b}.
\end{itemize}

\paragraph{Step 2b.}

\begin{itemize}
\item Set $\tilde{W}_n := W_n - \varepsilon \nabla\mathcal{J}(W_n)$, where
the gradient descent step size $\varepsilon$ is chosen numerically.

\item Set $\tilde{V}_n: = \tilde{W}_n + F$ and compute the corresponding
scalar function $\tilde{v}_n$.

\item Compute $a_n(x)$ from $\tilde{v}_n$ using the real part of (\ref{hev}) with $%
k=\underline k$.

\item Compute $u(x,k)$ by solving the direct problem~\eqref{Helm}--%
\eqref{radiation} with $a(x) := a_n(x)$, and set $u_{n+1}(x,k):=u(x,k)$.

\item Compute $v_{n+1}(x,k)$ from $u_{n+1}(x,k)$ using (\ref{p}) and (\ref{v}%
), and then compute $V_{n+1}$.

\item Set $n := n+1$ and return to \textbf{Step 2a}.
\end{itemize}


\paragraph{Step 3.}

Set $V := V_n$, compute $v(x,k)$ from $V$ and compute $a(x)$ using the real part of 
(\ref{hev}) with $k = \underline{k}$. 

\begin{remark}
We observe from the numerical performance of the algorithm that solving the
direct problem to update $V_{n}$ helps the cost functional decrease faster
with respect to iterations. This is important to the algorithm since the
cost functional decreases very slowly after the first iteration without this
update.
\end{remark}

\section{Numerical study}

In this section, we describe some important details of the numerical
implementation of the above algorithm and present some numerical
reconstruction results. The first step of the algorithm is to construct the
data carrier $F$. Recall that our computational domain $\Omega =(-R,R)^{2}$.
For $0<\xi <R$, define the $\chi _{0}(t),$ 
\begin{equation*}
\chi _{0}(t)=\left\{ 
\begin{array}{ll}
\text{exp}\left( -\frac{R}{t+\xi }\right) , & t>-\xi  \\ 
0, & t\leq -\xi 
\end{array}%
\right. 
\end{equation*}%
and then set%
\begin{equation*}
\chi (t)=\frac{\chi _{0}(t)}{\chi _{0}(t)+\chi _{0}(R-t-2\xi )}.
\end{equation*}%
Then $\chi (t)=0$ for $-R<t\leq -\xi $, $\chi (t)=1$ for $R-\xi \leq t<R$,
and $\chi $ is a smooth transition from $0$ to $1$ on $[-\xi ,R-\xi ]$. Set 
\begin{equation*}
f(x,k)=\left[ \tilde{g}_{0}(x,k)+(x_{2}-R)\tilde{g}_{1}(x,k)\right] \chi
(x_{2}),
\end{equation*}%
where  $ \tilde{g}_{0}$ and $ \tilde{g}_1$ are the Cauchy data for $v(x,k)$,  which means   $v = \tilde{g}_{0}$ on $\Gamma$ and $\partial_{x_2} v = \tilde{g}_{1}$ on $\Gamma$. Obviously, these data can be computed from the
data $g_{0}$ and $g_{1}$ in (\ref{data2})--(\ref{data3}) using the relation between $u$ and $v$ in
\eqref{p} and~\eqref{v}. 
Then the function  $f$ satisfies the boundary conditions $f=\tilde{g}_{0}$
on $\Gamma $ and $\partial _{x_{2}}f=\tilde{g}_{1}$ on $\Gamma $. Thus, the
corresponding vector valued function $F(x)=[f_{1}(x)\ f_{2}(x)\ \dots \
f_{N}(x)]^{T}$ containing the Fourier coefficients of $f$ with respect to
that truncated Fourier basis satisfies the boundary conditions (\ref{dd})--(%
\ref{nd}). Actually by the definition of $\chi $ the function $f(x,k)$ is
zero in $(-R,-\xi ]$, which means that we mainly seek for the scattering
object in the upper part $(-R,R)\times \lbrack -\xi ,R)$ of the square $%
\Omega $ since the Cauchy data are given on the top boundary $\Gamma $. 
Indeed, since the Cauchy problem (\ref{eqn})--(\ref{nd}) is
unstable, then it is unlikely that even after the regularization, which we
do here, one could image scattering objects well, if they are located far
from the measurement side $\Gamma $. For the parameter $\xi $ we
choose $\xi =R/10$ in the numerical implementation.

For the implementation of the algorithm, we first discretize the
computational domain $\Omega $ into $(N_{x}+1)\times (N_{x}+1)$ uniform grid
points $x_{ij}=(x_{j},y_{i}),1\leq i,j\leq N_{x}+1$, where the mesh size is $%
h_{x}$. The wave number interval $[\underline{k},\overline{k}]$ is divided
into $N_{k}$ uniform subintervals, where $k_{1},k_{2},\dots ,k_{N_{k}}$ are
the midpoints and $h_{k}$ is the length of each subinterval. Define the
lined up index as follows 
\begin{equation*}
\mathfrak{m}=\mathfrak{m}(i,j,r)=i+(j-1)(N_{x}+1)+(r-1)(N_{x}+1)^{2},\quad
1\leq i,j\leq N_{x}+1,\ 1\leq r\leq N.
\end{equation*}%
In this section, using the lined up index $\mathfrak{m,}$ we write vector
valued functions at grid points $x_{ij}$ as a column vector without changing
notations. For instance, for $U(x)=[u_{1}(x)\ u_{2}(x)\ \dots \ u_{N}(x)]^{T}
$, we have 
\begin{equation*}
U=[u_{\mathfrak{m}}],\quad 1\leq \mathfrak{m}\leq (N_{x}+1)^{2}N,
\end{equation*}%
where 
\begin{equation*}
u_{\mathfrak{m}}=u_{\mathfrak{m}(i,j,r)}=u_{r}(x_{ij}).
\end{equation*}


Let $\hat{W}=W+F$ and set $W=[w_{\mathfrak{m}}]$ and $\hat{W}=[\hat{w}_{%
\mathfrak{m}}]$. The weighted Tikhonov-like functional $\mathcal{J}$ in~%
\eqref{cost} is discretized using finite differences as 
\begin{align*}
\mathcal{J}(W)
&=h_{x}^{2}\sum_{m=1}^{N}\sum_{j=2}^{N_{x}}\sum_{i=2}^{N_{x}}\left\vert
\sum_{r=1}^{N}\left[ \frac{d_{mr}}{h_{x}^{2}}\left( \hat{w}_{\mathfrak{m}%
(i+1,j,r)}+\hat{w}_{\mathfrak{m}(i-1,j,r)}+\hat{w}_{\mathfrak{m}(i,j+1,r)}+%
\hat{w}_{\mathfrak{m}(i,j-1,r)}-4\hat{w}_{\mathfrak{m}(i,j,r)}\right)
\right. \right.  \\
& \quad +\sum_{s=1}^{N}\frac{b_{mr}^{s}}{h_{x}^{2}}\left( \hat{w}_{\mathfrak{%
m}(i,j+1,r)}-\hat{w}_{\mathfrak{m}(i,j,r)}\right) \left( \hat{w}_{\mathfrak{m%
}(i,j+1,s)}-\hat{w}_{\mathfrak{m}(i,j,s)}\right)  \\
& \quad +\sum_{s=1}^{N}\frac{b_{mr}^{s}}{h_{x}^{2}}\left( \hat{w}_{\mathfrak{%
m}(i+1,j,r)}-\hat{w}_{\mathfrak{m}(i,j,r)}\right) \left( \hat{w}_{\mathfrak{m%
}(i+1,j,s)}-\hat{w}_{\mathfrak{m}(i,j,s)}\right)  \\
& \quad +\left. \left. \frac{s_{mr}}{h_{x}}\left( \hat{w}_{\mathfrak{m}%
(i+1,j,r)}-\hat{w}_{\mathfrak{m}(i,j,r)}\right) \right] \varphi
(x_{ij})\right\vert ^{2} \\
& +\rho
h_{x}^{2}\sum_{m=1}^{N}\sum_{j=1}^{N_{x}+1}\sum_{i=1}^{N_{x}+1}\left\vert w_{%
\mathfrak{m}(i,j,m)}\right\vert ^{2} \\
& +\rho h_{x}^{2}\sum_{m=1}^{N}\sum_{j=2}^{N_{x}}\sum_{i=2}^{N_{x}}\left[
\left\vert \frac{w_{\mathfrak{m}(i,j+1,m)}-w_{\mathfrak{m}(i,j,m)}}{h_{x}}%
\right\vert ^{2}+\left\vert \frac{w_{\mathfrak{m}(i+1,j,m)}-w_{\mathfrak{m}%
(i,j,m)}}{h_{x}}\right\vert ^{2}\right.  \\
& \quad +\left\vert \frac{w_{\mathfrak{m}(i,j+1,m)}-2w_{\mathfrak{m}%
(i,j,m)}+w_{\mathfrak{m}(i,j-1,m)}}{h_{x}^{2}}\right\vert ^{2}+\left\vert 
\frac{w_{\mathfrak{m}(i+1,j,m)}-2w_{\mathfrak{m}(i,j,m)}+w_{\mathfrak{m}%
(i-1,j,m)}}{h_{x}^{2}}\right\vert ^{2} \\
& \quad +\left. 2\left\vert \frac{w_{\mathfrak{m}(i+1,j+1,m)}-w_{\mathfrak{m}%
(i-1,j+1,m)}-w_{\mathfrak{m}(i+1,j-1,m)}+w_{\mathfrak{m}(i-1,j-1,m)}}{%
h_{x}^{2}}\right\vert ^{2}\right]  \\
& +\alpha _{1}h_{x}\sum_{m=1}^{N}\sum_{j=1}^{N_{x}+1}\left\vert w_{\mathfrak{%
m}(N_{x}+1,j,m)}\right\vert ^{2}+\alpha
_{2}h_{x}\sum_{m=1}^{N}\sum_{j=2}^{N_{x}}\left\vert \frac{w_{\mathfrak{m}%
(N_{x}+1,j,m)}-w_{\mathfrak{m}(N_{x},j,m)}}{h_{x}}\right\vert ^{2}.
\end{align*}%
We now describe how to compute $\nabla \mathcal{J}(W)$ for complex valued
vector function $W$. Recall that if $z$ is a complex variable, $\mathbf{z}%
=(z_{1},z_{2},\dots ,z_{M})$ is a complex vector and $h(\mathbf{z})$ is a
complex valued function then we have (see~\cite{Kreut2009}) 
\begin{align*}
& i.\quad \frac{\partial }{\partial z}|z|^{2}=\frac{\partial }{\partial z}(z%
\bar{z})=\bar{z}, \\
& ii.\quad \frac{\partial h}{\partial \mathbf{z}}(\mathbf{z})=\left[ \frac{%
\partial h}{\partial z_{1}}\ \frac{\partial h}{\partial z_{2}}\ \dots \ 
\frac{\partial h}{\partial z_{M}}\right] , \\
& iii.\quad \nabla h(\mathbf{z})=\left( \overline{\frac{\partial h}{\partial 
\mathbf{z}}(\mathbf{z})}\right) ^{T}.
\end{align*}%
Thus, by the chain rule we have 
\begin{equation*}
\frac{\partial }{\partial \mathbf{z}}|h|^{2}(\mathbf{z})=\overline{h(\mathbf{%
z})}\ \frac{\partial h}{\partial \mathbf{z}}(\mathbf{z}).
\end{equation*}%
Treating $\mathcal{J}(W)$ as a function of $(N_{x}+1)^{2}N$ complex
variables and applying all of the above to each of its summands, we are able
to compute $\frac{\partial \mathcal{J}}{\partial W}(W)$ and thus obtain $%
\nabla \mathcal{J}(W)$ using ($iii$).

We need a numerical solver for the direct problem~\eqref{Helm}--%
\eqref{radiation} in Step 2b of the reconstruction algorithm and to generate
synthetic scattering data for the numerical study. It is well known that the
direct problem~\eqref{Helm}--\eqref{radiation} is equivalent to the
Lippmann-Schwinger integral equation 
\begin{equation}
u(x,k)=u_{\mathrm{in}}(x,k)+k^{2}\int_{\Omega }\frac{i}{4}%
H_{0}^{(1)}(k\lvert x-y\rvert )a(y)u(y,k)dy,  \label{ls}
\end{equation}%
where $H_{0}^{(1)}$ is the Hankel function of the first kind of order 0, see~%
\cite{Colto2013}. We exploit the numerical method studied in~\cite{Saran2002}
to solve this integral equation to generate the Cauchy data $g_{0}(x,k)$ and 
$g_{1}(x,k)$ for the inverse problem. Note that the numerical method studied
in~\cite{Saran2002} assumes smooth coefficients. Its extension to the case
of discontinuous coefficients is studied in~\cite{Lechl2014} which can be
adapted to our discontinuous coefficient examples in this section. We also
add an artificial random noise to the  data 
\begin{equation*}
g_{j}(x,k)=g_{j}(x,k)+\delta \lVert g_{j}\rVert _{L^{2}}\mathcal{N}%
_{j}(x,k),\quad j=0,1,
\end{equation*}%
where $\delta $ is the noise level and $\mathcal{N}_{j}$ are functions
taking random complex values and satisfy $\lVert \mathcal{N}_{j}\rVert
_{L^{2}}=1$. We consider $5\%$ noise in the Cauchy backscatter data which
means $\delta =0.05$ in our numerical examples.

In all numerical examples presented in this section the computational domain
is chosen as $\Omega =(-0.8,0.8)^{2}$, where $N_{x}=28$. This means $\Omega $
is uniformly discretized into $29^{2}$ points. The interval of wave numbers
is $k\in \lbrack \underline{k},\overline{k}]=[0.5,2]$, where $N_{k}=50$. We
have found that $N=4$ in (\ref{100})\textbf{\ }is sufficient for the Fourier
series truncation, also, see~\cite{Khoa2019} for a similar choice. We
generate the multifrequency data for the incident plane wave 
\begin{equation*}
u_{\mathrm{in}}(x)=e^{-ikx_{2}},\quad k\in \lbrack 0.5,2].
\end{equation*}%
In the Carleman Weight Function $\varphi (x)$ (\ref{cw}), we choose $\lambda
=5,s=1.$ This means that 
\begin{equation*}
\varphi (x)=e^{-5(x_{2}-1)^{2}}.
\end{equation*}%
This choice was made by trial and error, so as choices of all other
parameters used in this section. We refer to works on the convexification~\cite{Kliba2018a, Kliba2019a, Kliba2017, Kliba2019c, Kliba2019d, Khoa2019}
 for choices of smaller $\lambda \in \left[ 1,3\right]
$. As to our choice of $\lambda ,$ it seems to give us the optimal results
for our numerical examples of this section. The step size $\varepsilon $ of
the gradient descent method and the regularization parameters $\rho ,\alpha
_{1},\alpha _{2}$ were chosen as: 
\begin{equation*}
\varepsilon =\alpha _{1}=10^{-3},\quad \rho =\alpha _{2}=10^{-5}.
\end{equation*}%
The tolerance number of the iterations is set up to be $10^{-3}$. It happens
in all our numerical examples that the algorithm stops within 10 iterations.
One can observe in the numerical examples that the value of the minimized
functional $\mathcal{J}(W)$ does not change much after 8 or 9 iterations.
Since we are interested in $a(x)\geq 0$, in the final iteration we assign $%
a(x)$ to be zero in the area in which it takes negative values. By our
numerical experience, this area is typically below the reconstructed
scatterer.

It is important to mention that the initial guess $W_{0}\equiv 0$
 for $W$ in all numerical examples below. This goes along
well with our theory (to be published) which guarantees that our algorithm
converges to the correct solution starting from any point of the ball $%
B( M) $ defined in (\ref{201}), see item 3 in the end of
Introduction. This certainly a significant advantage of our method, compared
with locally convergent optimization approaches, which typically need a
strong \emph{a priori} knowledge of the scatterer. Such a knowledge, however, is
rarely available in applications.

%

\subsection{Numerical example 1}

\begin{figure}[h!]
\centering
\subfloat[True $a(x)$]{\includegraphics[width=6cm]{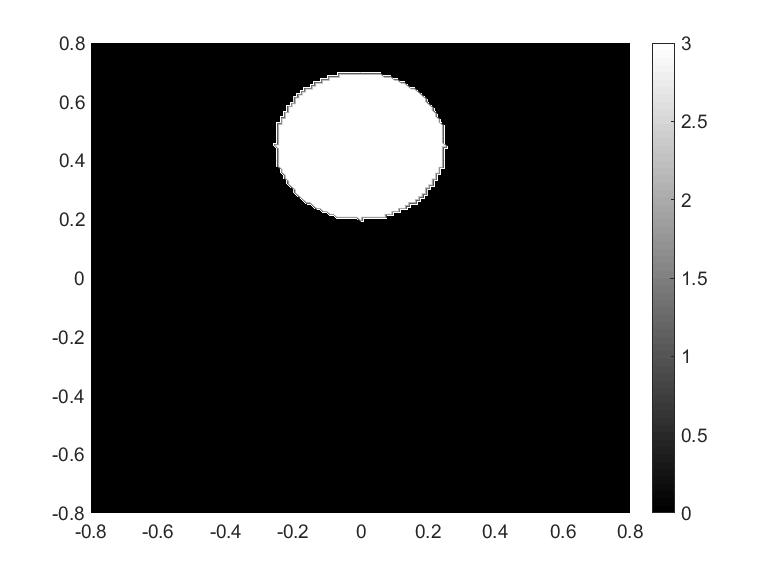}}\hspace{%
-0.0cm} 
\subfloat[The cost functional]{\includegraphics[width=6cm]{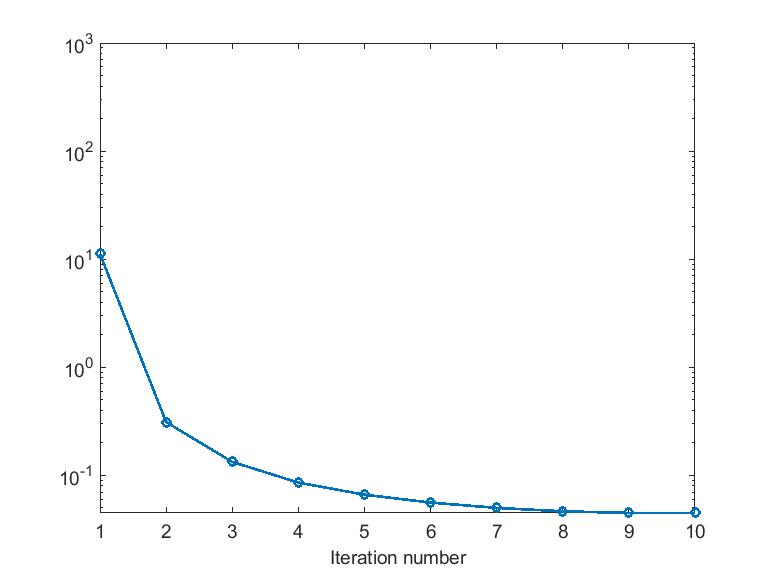}}\hspace{-0.0cm} 
\subfloat[Reconstruction]{\includegraphics[width=6cm]{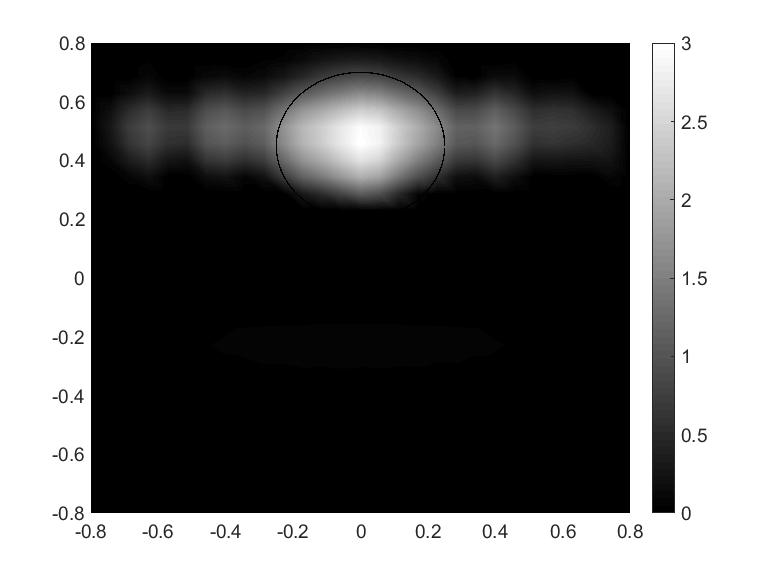}} 
\hspace{-0.0cm} 
\subfloat[View at $x_2 = 0.45$]{\includegraphics[width=6cm]{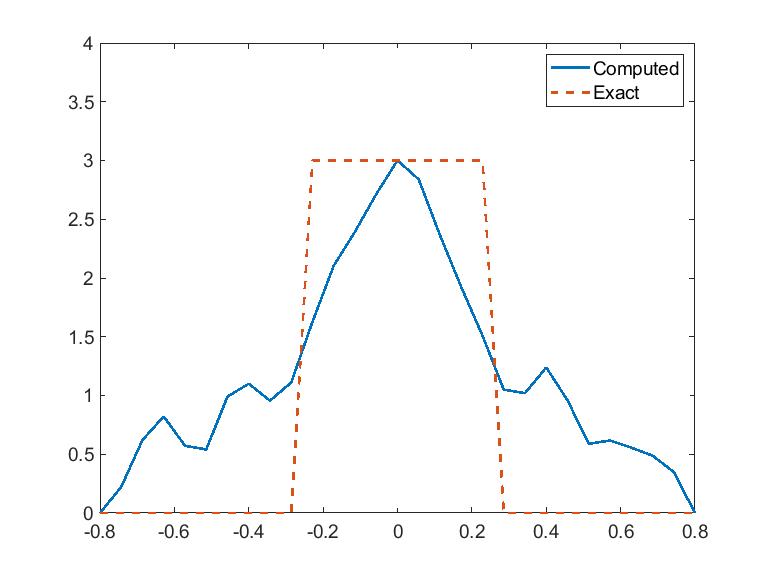}} 
\caption{Reconstruction of one scattering disk characterized by $a(x) = 3$.}
\label{fi1}
\end{figure}

In this example we consider a single scattering disk characterized by the
coefficient $a(x)$ which equals 3 inside the disk and zero elsewhere. We can
see from the reconstruction result in Figure~\ref{fi1} that the location and
and the maximal value of $a(x)$ are well reconstructed. It seems to us that
the shape of the scattering object is not well-reconstructed because the
backscatter data are generated by incident plane waves with a fixed
direction, also, see~\cite{Nguye2017,Nguye2018, Kliba2019a, Kliba2019d} for
similar results. Convergence of the algorithm can be observed from Figure~%
\ref{fi1}(b). From our numerical experience the cost functional $\mathcal{J}$
does not decrease much after 8 or 9 iterations, see also Figure~\ref{fi1}(b).

\begin{figure}[h!]
\centering
\subfloat[Reconstruction result at the 6th
iteration]{\includegraphics[width=6cm]{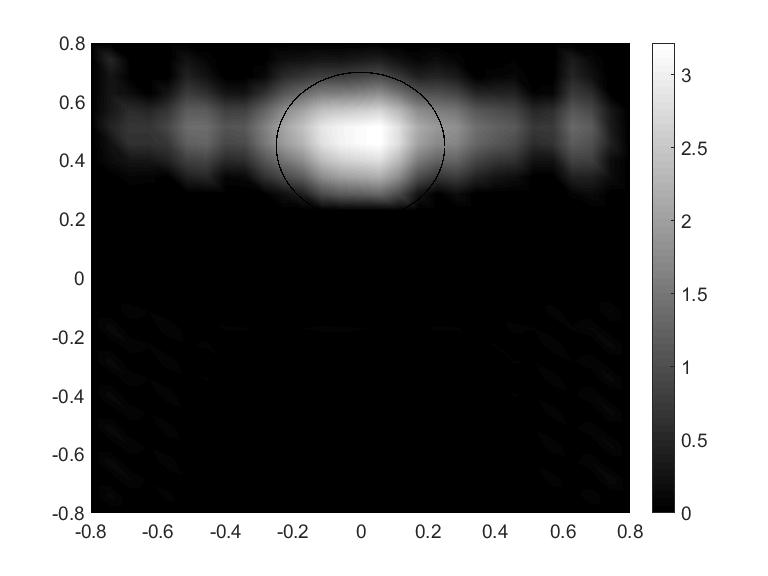}}\hspace{-0.0cm} %
\subfloat[The cost functional]{\includegraphics[width=6cm]{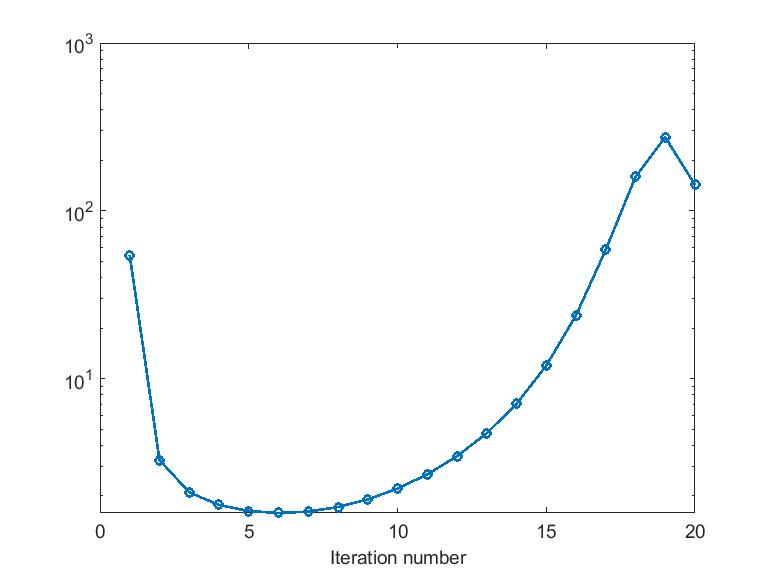}}%
\hspace{-0.0cm} 
\caption{Reconstruction of the scattering disk in Figure~\protect\ref{fi1}%
(a) in which the cost functional does not involve the Carleman Weight
Function.}
\label{fi0}
\end{figure}

Now with the numerical result in Figure~\ref{fi0} we want to indicate the
importance of the Carleman Weight Function for our numerical algorithm. The
algorithm does not converge when the cost functional $\mathcal{J}$ does not
involve the Carleman weight function. Firstly, the error between the cost
functionals at two consecutive iterations is never smaller than the
tolerance number $10^{-3}$ like what we have when the Carleman Weight
Function is present. Therefore, the iterations do not stop with the chosen
tolerance. Secondly, the cost functional $\mathcal{J}$ starts to increase
after a certain number of iterations. Figure~\ref{fi0}(a) presents the
reconstruction result at the sixth iteration where the cost functional
obtains its smallest value among 20 iterations. However, this result is not
as good as that of Figure~\ref{fi1}(c) where the Carleman Weight Function is
involved. Indeed, the artifact in Figure~\ref{fi0}(a) is slightly stronger and the reconstructed maximal value 
is 3.2157 while the maximal value of the reconstruction  in Figure~\ref{fi1}(c) is  3.0014.  
 Also for the next examples, the reconstruction results are always
better with the presence of the Carleman Weight Function in the cost
functional.

%

\subsection{Numerical example 2}

\begin{figure}[]
\centering
\subfloat[True $a(x)$]{\includegraphics[width=6cm]{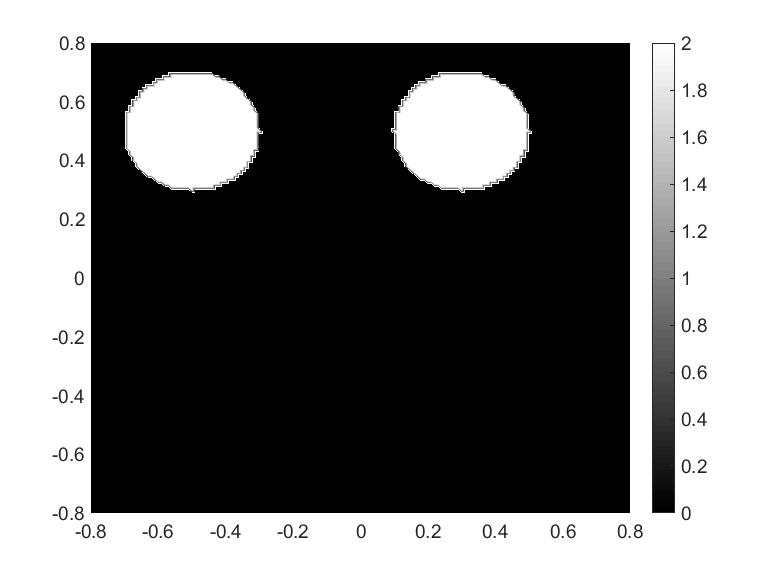}} \hspace{%
-0.0cm} 
\subfloat[The cost
functional]{\includegraphics[width=6cm]{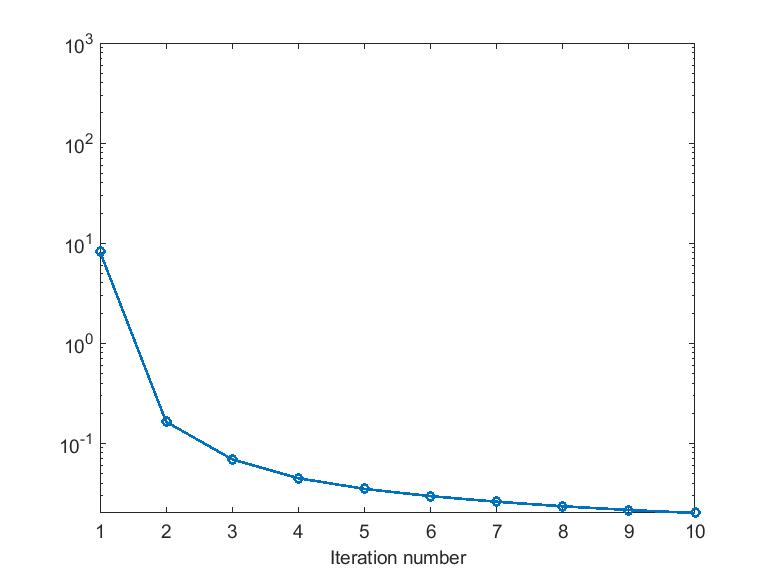}}\hspace{-0.0cm} 
\subfloat[Reconstruction]{\includegraphics[width=6cm]{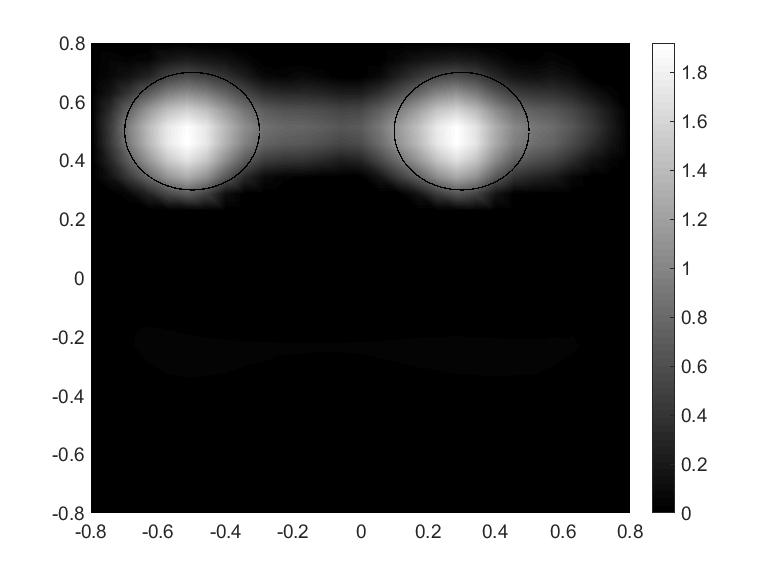}} 
\hspace{-0.0cm} 
\subfloat[View at $x_2 =
0.45$.]{\includegraphics[width=6cm]{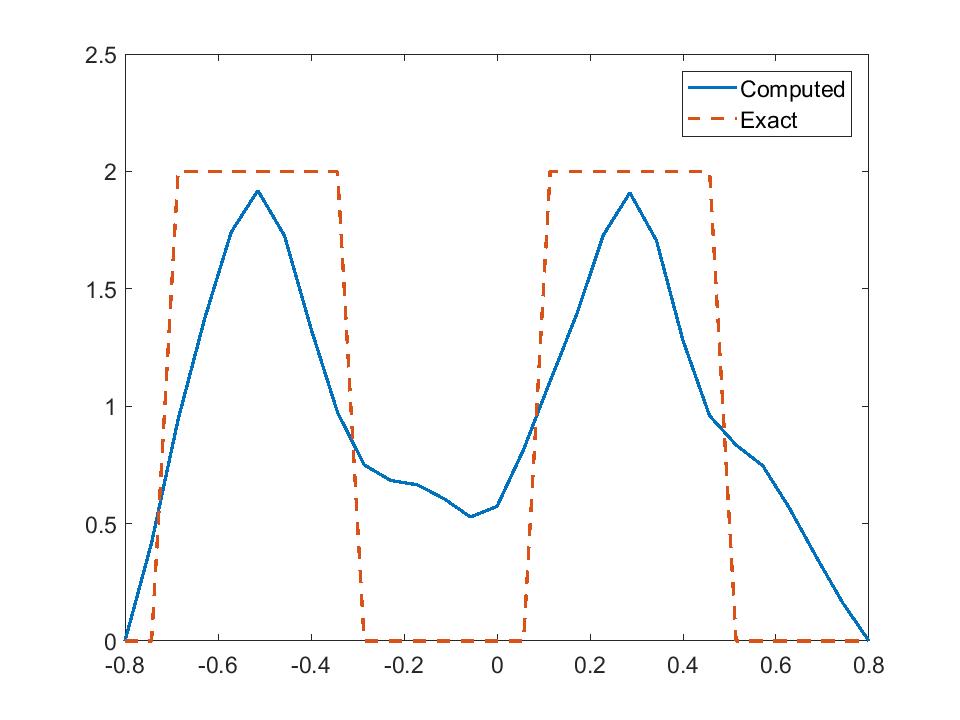}} 
\caption{Reconstruction of two similar scattering disks characterized by the
coefficient $a(x) = 2$. }
\label{fi2}
\end{figure}

\begin{figure}[]
\centering
\subfloat[True $a(x)$]{\includegraphics[width=6cm]{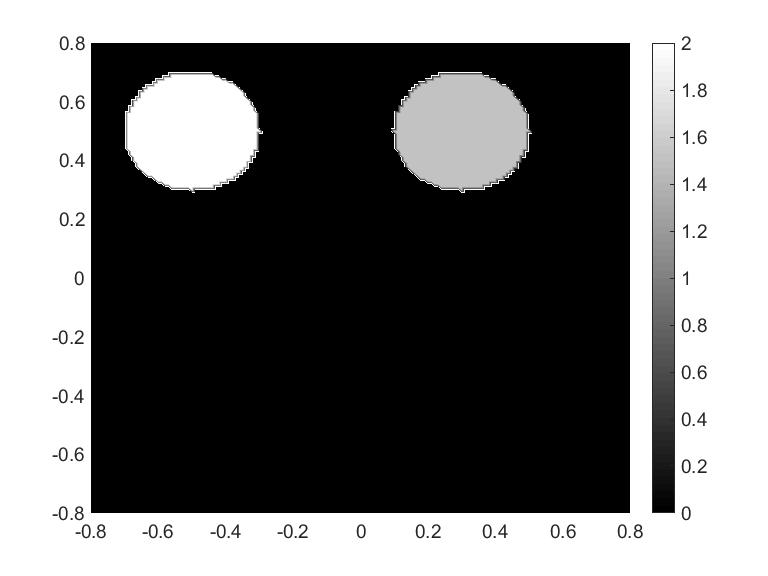}}\hspace{%
-0.0cm} 
\subfloat[The cost
functional]{\includegraphics[width=6cm]{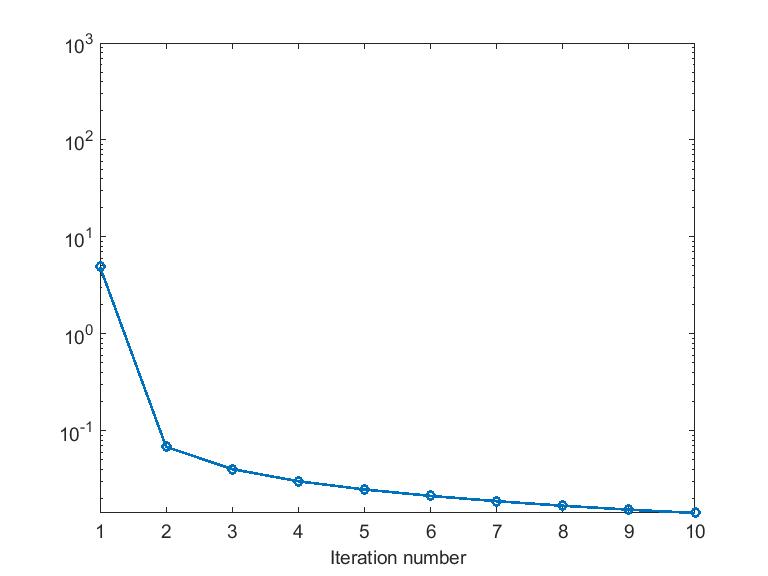}}\hspace{-0.0cm} 
\subfloat[Reconstruction]{\includegraphics[width=6cm]{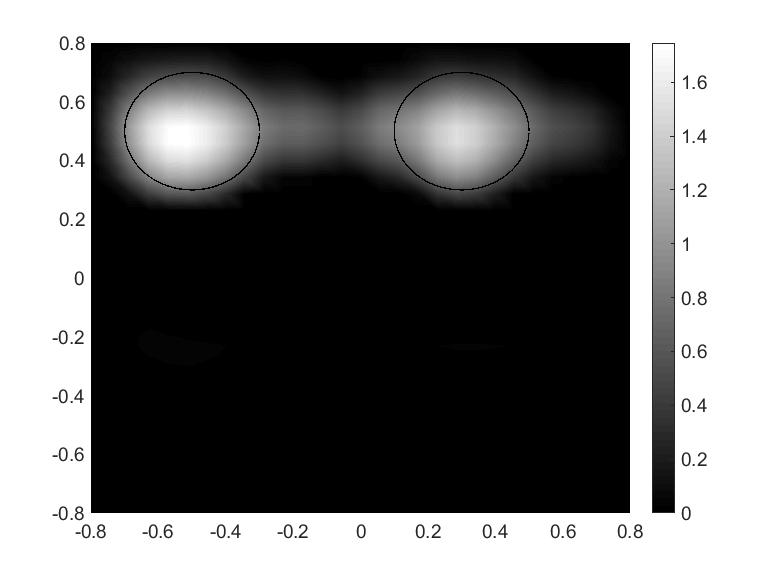}}%
\hspace{-0.0cm} 
\subfloat[View at $x_2 =
0.45$]{\includegraphics[width=6cm]{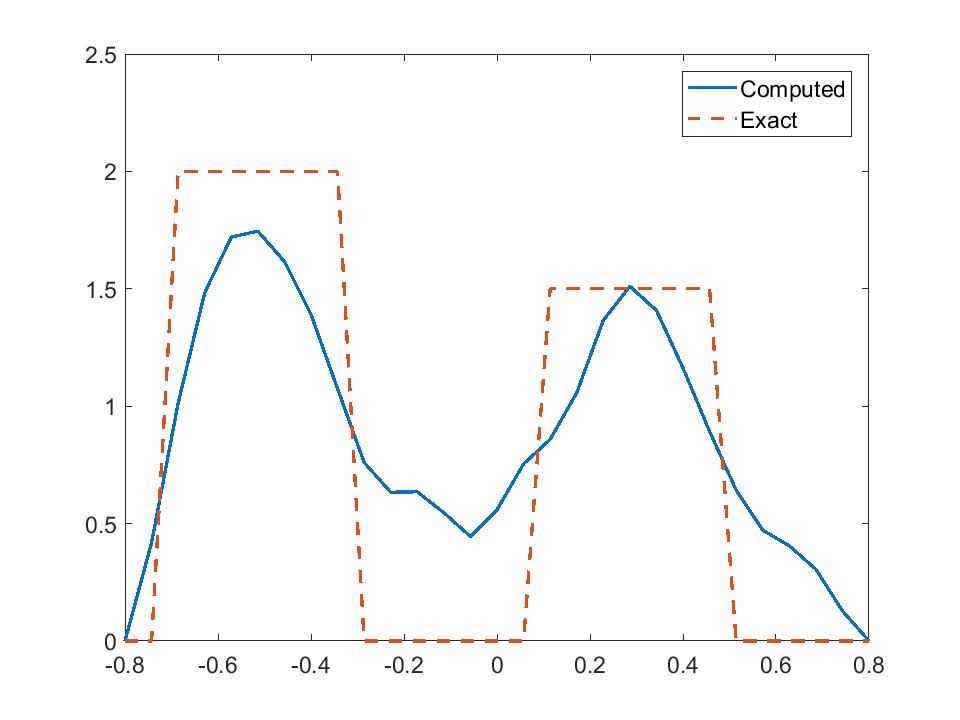}} 
\caption{Reconstruction of two scattering disks with different values. The
coefficient $a(x) = 2$ in the left scattering disk and $a(x) = 1.5$ in the
right scattering disk. }
\label{fi3}
\end{figure}

In this example we consider the case of two scattering disks. In the first
case in Figure~\ref{fi2}(a) two similar scattering disks are considered. The
reconstruction result in Figure~\ref{fi2}(c) again shows that the algorithm
is able to reconstruct very well the location and the maximal values of $a(x)
$ in this case. The cost functional decreases well within ten iterations,
see Figure~\ref{fi2}(b). The case of Figure~\ref{fi3}(a) is more challenging
since the maximal values of $a(x)$ in each scattering disk are different.
However, the algorithm can provide reasonable reconstruction results in
Figures~\ref{fi3}(c). One can clearly see the locations of the scattering
disks as well as two different maximal values of $a(x)$ on each disk. We
point out the algorithm in this case can reconstruct the scatterer
consisting of two components without using any \emph{\emph{a priori}} knowledge
about the number of components of the scatterer. 

\subsection{Numerical example 3}

\begin{figure}[]
\centering
\subfloat[True $a(x)$]{\includegraphics[width=6cm]{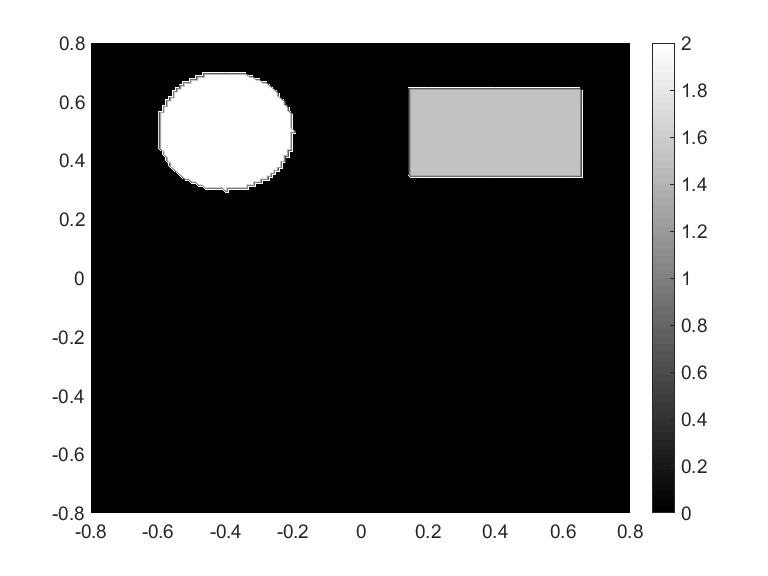}}\hspace{%
-0.0cm} 
\subfloat[The cost
functional]{\includegraphics[width=6cm]{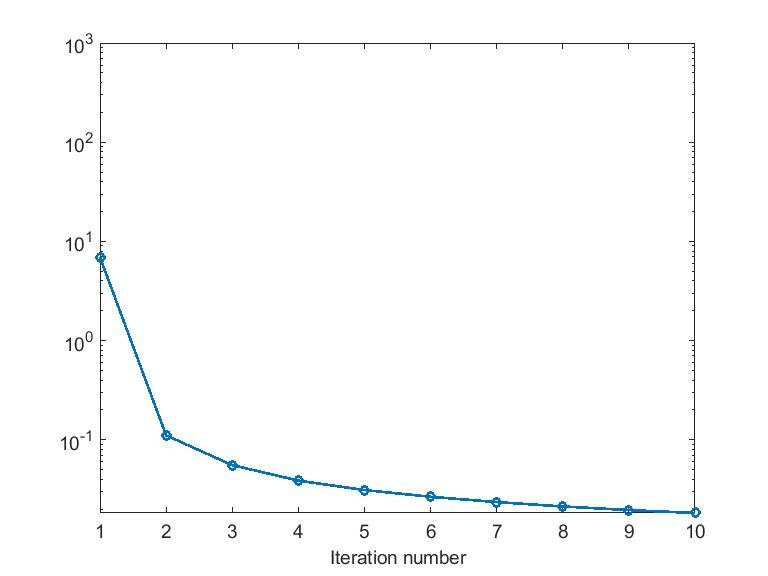}}\hspace{-0.0cm} 
\subfloat[Reconstruction]{\includegraphics[width=6cm]{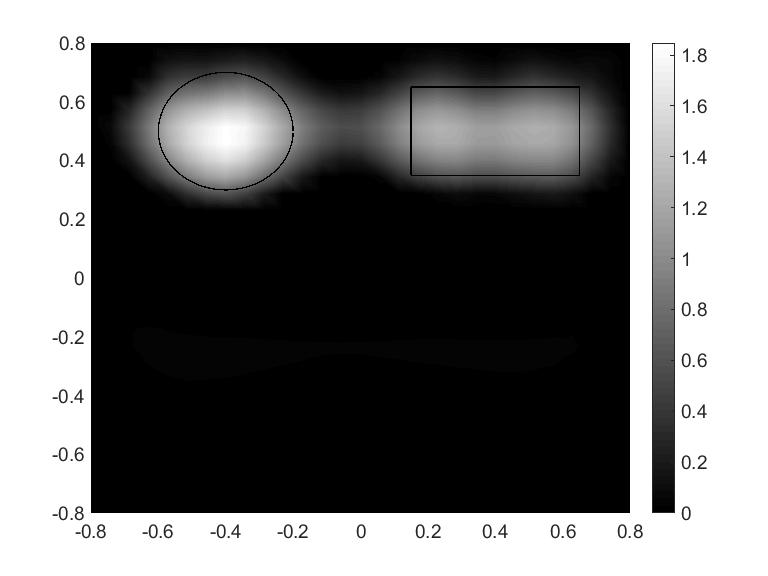}}%
\hspace{-0.0cm} 
\subfloat[View at $x_2 =
0.45$]{\includegraphics[width=6cm]{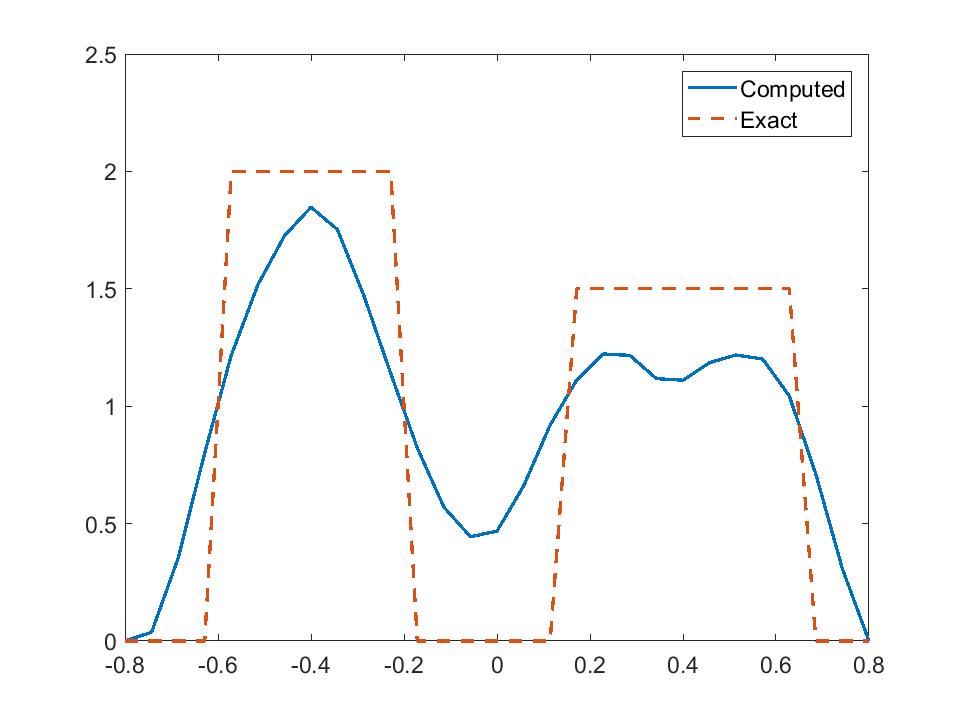}} 
\caption{Reconstruction of scattering objects with different shapes and
values. The coefficient $a(x) = 2$ in the left scattering object and $a(x) =
1.5$ in the right scattering object. }
\label{fi4}
\end{figure}

\begin{figure}[]
\centering
\subfloat[True $a(x)$]{\includegraphics[width=6cm]{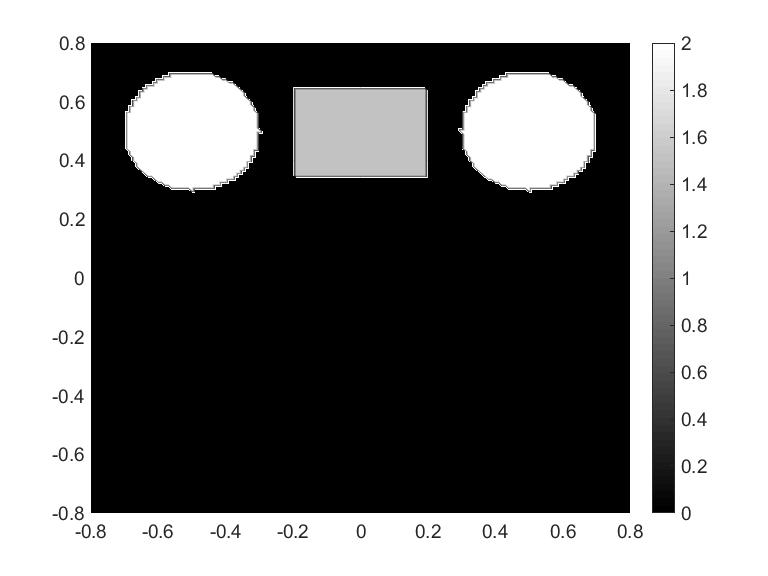}}\hspace{%
-0.0cm} 
\subfloat[The cost
functional]{\includegraphics[width=6cm]{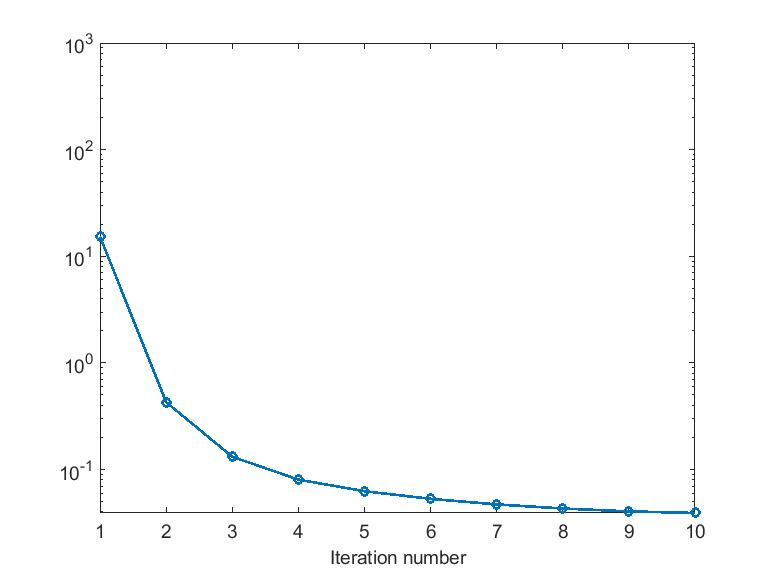}}\hspace{-0.0cm} 
\subfloat[Reconstruction]{\includegraphics[width=6cm]{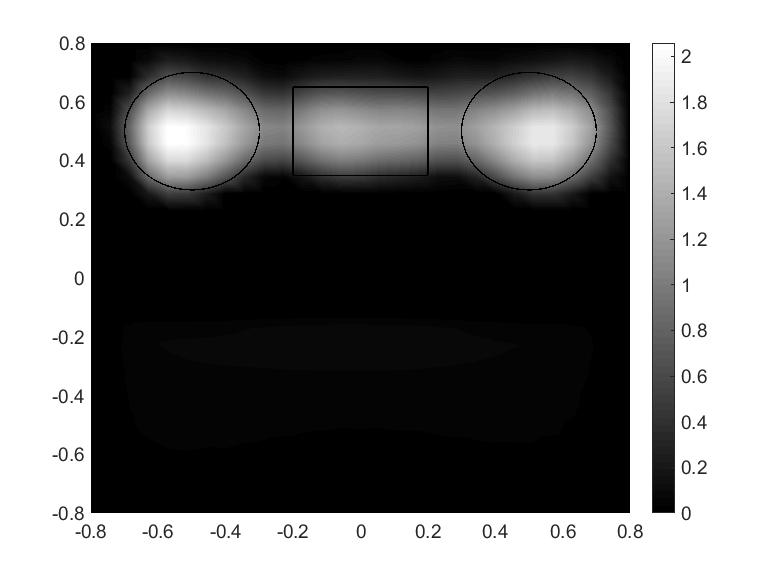}}%
\hspace{-0.0cm} 
\subfloat[View at $x_2 =
0.45$]{\includegraphics[width=6cm]{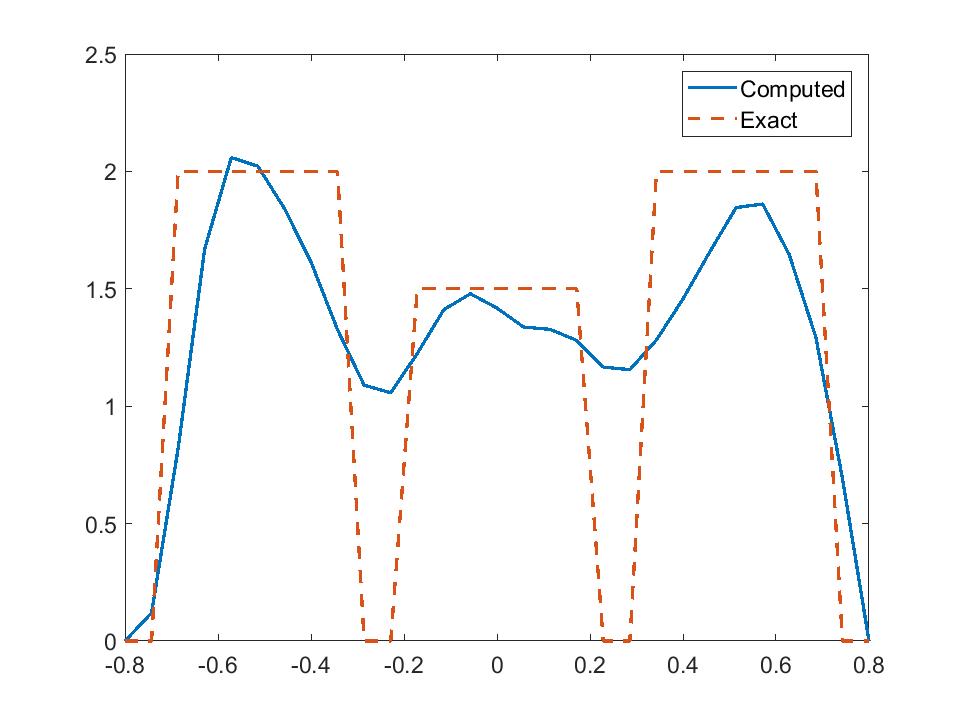}} 
\caption{Reconstruction of three scattering objects with different shapes
and values. The coefficient $a(x) = 2$ in the scattering disks and $a(x) =
1.5$ in the scattering rectangle. }
\label{fi5}
\end{figure}

In this example we consider the case of the coefficient $a(x)$ which has
different values in scattering objects of different shapes. This case is
thus more challenging than those of the first two examples. The scatterer in
Figure~\ref{fi4}(a) consists of a scattering disk in which $a(x)=2$ and a
scattering rectangle in which $a(x)=1.5$. The reconstruction result in
Figure~\ref{fi4}(c) shows that the algorithm again can compute the location
of the scatterer and the maximal values of $a(x)$ in each scattering object.
Particularly, we can also see pretty well a difference between the shape of
the disk and the rectangle in the reconstruction. The scatterer in Figure~%
\ref{fi5}(a) consists of two scattering disks in which $a(x)=2$ and a
scattering rectangle in which $a(x)=1.5$. The reconstruction result in
Figure~\ref{fi5}(c) provides the location and maximal values of the
scattering objects. However, the resolution of the reconstruction in this
case is not as good as that of the case of two objects since three
scattering objects are placed quite close to each other.

\section{Summary}
We have proposed a new version of the convexification numerical
reconstruction method for solving the coefficient inverse scattering problem
with multifrequency backscatter data associated with a single direction of the
incident plane wave. This method relies on an approximate reformulation of
the problem as the Cauchy problem for a system of coupled quasilinear
elliptic PDEs. The main ingredients for deriving this formulation are the
elimination of the coefficient from Helmholtz equation and the use of
truncated Fourier expansion for the total field. To solve the quasilinear
elliptic PDE system, we use a weighted quasi-reversibility method in which a
Carleman Weight Function is included in the weighted Tikhonov-like
functional. The numerical results show that our method is able to
efficiently compute the solution without using any \emph{a priori} information
about it. We have shown that values of the dielectric constants of
scatterers as well as locations of scatterers can be well reconstructed.
However, shapes are not reconstructed well. On the other hand, the recent
publication~\cite{Khoa2019} shows that all three components of
scatterers can be accurately reconstructed in the case when the point source
moves along a straight line and frequency is fixed.  

Overall, the main advantage of our algorithm, so as other above
cited versions of the convexification method, is its rigorously guaranteed
global convergence, as opposed to the local convergence of the conventional
optimization methods, see Introduction for the definition of the global
convergence. Theoretical analysis of the algorithm as well as its
three-dimensional extension will be addressed in forthcoming publications. 
That theoretical analysis will consists in detailed proofs of results
announced in items 1-3 of Introduction.\\

\noindent \textbf{Acknowledgement.} The authors are grateful to Loc Nguyen
and Khoa Vo from the University of North Carolina at Charlotte for helpful
discussions. The first and second authors are partially supported by NSF
grant DMS-1812693. The third author is supported by US Army Research
Laboratory and US Army Research Office grant W911NF-19-1-0044.

\bibliographystyle{elsart-num-sort.bst}
\bibliography{ip-biblio}

\providecommand{\noopsort}[1]{}
\begin{thebibliography}{10}
\expandafter\ifx\csname url\endcsname\relax
  \def\url#1{\texttt{#1}}\fi
\expandafter\ifx\csname urlprefix\endcsname\relax\def\urlprefix{URL }\fi

\bibitem{Bak}
A.~B. Bakushinsky, M.~V. Klibanov, N.~A. Koshev, Carleman weight functions for
  a globally convergent numerical method for ill-posed {C}auchy problems for
  some quasilinear {PDE}s, Nonlinear Analysis: Real World Applications 34
  (2017) 201--224.

\bibitem{Bakus2004}
A.~B. Bakushinsky, M.~Y. Kokurin, Iterative Methods for Approximate Solutions
  of Inverse Problems, Springer Verlag, 2004.

\bibitem{Baud1}
L.~Baudouin, M.~de~Buhan, S.~Ervedoza, Convergent algorithm based on carleman
  estimates for the recovery of a potential in the wave equation, SIAM J.
  Numer. Anal. 55 (2017) 1578--1613.

\bibitem{Baud2}
L.~Baudouin, M.~de~Buhan, S.~Ervedoza, A.~Osses, Carleman-based reconstruction
  algorithm for the waves, preprint,
  https://hal.archives-ouvertes.fr/hal-02458787.

\bibitem{Beili2012}
L.~Beilina, M.~V. Klibanov, Approximate Global Convergence and Adaptivity for
  Coefficient Inverse Problems, Springer, New York, 2012.

\bibitem{BK}
L.~Beilina, M.~V. Klibanov, Globally strongly convex cost functional for a
  coefficient inverse problem, Nonlinear Anal. Real World Appl. 22 (2015)
  272--288.

\bibitem{BBS}
M.~Boulakia, M.~de~Buhan, E.~Schwindt, Numerical reconstruction based on
  carleman estimates of a source term in a reaction-diffusion equation,
  preprint, https://hal.archives-ouvertes.fr/hal-02185889.

\bibitem{Bukhg1981}
A.~L. Bukhgeim, M.~V. Klibanov, Global uniqueness of a class of
  multidimensional inverse problems, Soviet Math. Dokl. 24 (1981) 244--247.

\bibitem{Cakon2006}
F.~Cakoni, D.~Colton, Qualitative Methods in Inverse Scattering Theory. An
  Introduction., Springer, Berlin, 2006.

\bibitem{Chavent}
G.~Chavent, Nonlinear {L}east {S}quares for {I}nverse {P}roblems: {T}heoretical
  {F}oundations and {S}tep-by-{S}tep {G}uide for {A}pplications, {S}cientic
  {C}omputation, Springer, New York, 2009.

\bibitem{Colto2013}
D.~L. Colton, R.~Kress, Inverse acoustic and electromagnetic scattering theory,
  3rd ed., Springer, 2013.

\bibitem{Engl1996}
H.~W. Engl, M.~Hanke, A.~Neubauer, Regularization of inverse problems, Kluwer
  Acad. Publ., Dordrecht, Netherlands, 1996.

\bibitem{Gonch1}
A.~Goncharsky, S.~Romanov, A method of solving the coefficient inverse problems
  of wave tomography, Comput. Math. Appl. 77 (2019) 967--980.

\bibitem{Gonch2}
A.~Goncharsky, S.~Romanov, S.~Seryozhnikov, Low-frequency ultrasonic
  tomography: mathematical methods and experimental results, Moscow University
  Physics Bulletin 74 (2019) 43--51.

\bibitem{Khoa2019}
V.~Khoa, M.~V. Klibanov, L.~Nguyen, Convexification for a 3{D} inverse
  scattering problem with the moving point source, Submitted
  (arXiv:1911.10289).

\bibitem{Kirsc2008}
A.~Kirsch, N.~Grinberg, The Factorization Method for Inverse Problems, Oxford
  Lecture Series in Mathematics and its Applications 36, Oxford University
  Press, 2008.

\bibitem{Klib97}
M.~Klibanov, Global convexity in a three-dimensional inverse acoustic problem,
  SIAM J. Math. Anal. 28 (1997) 1371--1388.

\bibitem{Klib95}
M.~Klibanov, O.~Ioussoupova, Uniform strict convexity of a cost functional for
  three-dimensional inverse scattering problem, SIAM J. Appl. Math. 26 (1995)
  147--179.

\bibitem{Ksurvey}
M.~V. Klibanov, Carleman estimates for global uniqueness, stability and
  numerical methods for coefficient inverse problems, J. Inverse Ill-Posed
  Probl. 21 (2013) 477--560.

\bibitem{Klib2015}
M.~V. Klibanov, Carleman weight functions for solving ill-posed cauchy problems
  for quasilinear pdes, Inverse Problems 31 (2015) 125007.

\bibitem{Kliba2017}
M.~V. Klibanov, Convexification of restricted {D}irichlet-to-{N}eumann map, J.
  Inverse Ill-Posed Probl. 25 (2017) 669--685.

\bibitem{KK}
M.~V. Klibanov, V.~G. Kamburg, Globally strictly convex cost functional for an
  inverse parabolic problem, Math. Meth. Appl. Sci.. 39 (2015) 930--940.

\bibitem{Kliba2019d}
M.~V. Klibanov, A.~E. Kolesov, Convexification of a 3-{D} coefficient inverse
  scattering problem, Comput. Math. Appl. 77 (2019) 1681--1702.

\bibitem{Kliba2019a}
M.~V. Klibanov, A.~E. Kolesov, D.-L. Nguyen, Convexification method for an
  inverse scattering problem and its performance for experimental backscatter
  data for buried objects, SIAM J. Imaging Sci. 12 (2019) 576--603.

\bibitem{Kliba2018a}
M.~V. Klibanov, A.~E. Kolesov, A.~Sullivan, L.~Nguyen, A new version of the
  convexification method for a 1{D} coefficient inverse problem with
  experimental data, Inverse Problems 34 (2018) 115014.

\bibitem{Kliba2019c}
M.~V. Klibanov, J.~Li, W.~Zhang, Convexification for the inversion of a time
  dependent wave front in a heterogeneous medium, SIAM J. Appl. Math. 79 (2019)
  1722--1747.

\bibitem{Kliba2018}
M.~V. Klibanov, D.-L. Nguyen, L.~H. Nguyen, H.~Liu, A globally convergent
  numerical method for a 3{D} coefficient inverse problem with a single
  measurement of multi-frequency data, Inverse Probl. Imaging 12 (2018)
  493--523.

\bibitem{KT}
M.~V. Klibanov, A.~Timonov, Carleman Estimates for Coefficient Inverse Problems
  and Numerical Applications, VSP, Utrecht, 2004.

\bibitem{Koles2017}
A.~E. Kolesov, M.~V. Klibanov, L.~H. Nguyen, D.-L. Nguyen, N.~T. Th\`anh,
  Single measurement experimental data for an inverse medium problem inverted
  by a multi-frequency globally convergent numerical method, Appl. Numer. Math.
  120 (2017) 176--196.

\bibitem{Kreut2009}
K.~Kreutz-Delgado, The {C}omplex {G}radient {O}perator and the {CR} {C}alculus,
  Tech. rep., ArXiv:0906.4835v1 (2009).

\bibitem{LN}
T.~T. Le, L.~H. Nguyen, A convergent numerical method to recover the initial
  condition of nonlinear parabolic equations from lateral cauchy data,
  preprint, arXiv: 1910.05584.

\bibitem{Lechl2014}
A.~Lechleiter, D.-L. Nguyen, A trigonometric {G}alerkin method for volume
  integral equations arising in {TM} grating scattering, Adv. Comput. Math. 40
  (2014) 1--25.

\bibitem{Nguye2017}
D.-L. Nguyen, M.~V. Klibanov, L.~Nguyen, A.~E. Kolesov, M.~A. Fiddy, H.~Liu,
  Numerical solution for a coefficient inverse problem with multi-frequency
  experimental raw data by a globally convergent algorithm, J. Comput. Phys.
  345 (2017) 17--32.

\bibitem{Nguye2018}
D.-L. Nguyen, M.~V. Klibanov, L.~H. Nguyen, M.~Fiddy, Imaging of buried objects
  from multi-frequency experimental data using a globally convergent inversion
  method, J. Inverse Ill-Posed Probl. 26 (2018) 501--522.

\bibitem{Saran2002}
J.~Saranen, G.~Vainikko, Periodic integral and pseudodifferential equations
  with numerical approximation, Springer, 2002.

\end{thebibliography}

\end{document}